\documentclass{article}

\usepackage{amssymb} 
\usepackage{amsmath} 
\usepackage{epsfig}

\numberwithin{equation}{section}

\title{\bf Energy of tori of revolution}
\author{Hiroki Funaba and Jun O'Hara \\ 
\smallskip
{\small Department of Mathematics and Information Sciences, Tokyo Metropolitan University}\\ 
\smallskip
{\small E-mail:ohara@tmu.ac.jp }}

\begin{document}

\maketitle

\begin{abstract} 
We show that the surface energy introduced by Auckly and Sadun attains the minimum value at 
the Clifford torus among tori of revolution. 
%
\end{abstract}

\section{Introduction and main result}
%
\setlength\arraycolsep{1pt}
Recently Fernando C. Marques and Andr\'e Neves \cite{MN} proved the Willmore conjecture, namely, they showed the following. 
Let $\Sigma$ be an immersed torus in $\mathbb{R}^3$. Let $\kappa_1$ and $\kappa_2$ be principal curvatures. The Willmore functional is given by 
\[\mathcal{W}(\Sigma)=\int_\Sigma\left(\frac{\kappa_1+\kappa_2}2\right)^2\,d\Sigma=\int_\Sigma\left(\frac{\kappa_1-\kappa_2}2\right)^2\,d\Sigma,\]
where the second equation is the consequence of the Gauss-Bonnet theorem. 
It is known to be invariant under M\"obius transformations of $\mathbb{R}^3$. 
Then the Willmore conjecture, now the theorem of Marques and Neves, asserts that $\mathcal{W}(\Sigma)\ge2\pi^2$ and that the equality holds if and only if $\Sigma$ is a torus of revolution whose generating circle has radius $1$ and center at distance $\sqrt2$ from the axis of revolution up to a M\"obius transformation, in other words, if and only if $\Sigma$ is the image of a stereographic projection of the Clifford torus 
\[
\left\{(z_1,z_2)\in \mathbb C\times \mathbb C\,\Big|\,|z_1|=|z_2|=1/{\sqrt2}\,\right\}\subset
S^3=\left\{(z_1,z_2)\in \mathbb C\times \mathbb C\,\Big|\,|z_1|^2+|z_2|^2=1\,\right\}.\]

\smallskip

In this paper, we give another characterization of the Clifford torus using the surface energy introduced by David Auckly and Lorenzo Sadun (\cite{AS}), which is also invariant under M\"obius transformations. 
To be precise, we have not yet succeeded in proving that the Clifford torus gives the minimum energy among all the immersed tori\footnote{To show it, it suffice to show that the Clifford torus gives the minimum energy among all the embedded tori, since the energy blows up if a torus has a double point}. 
We only show that it gives the minimum energy among one-parameter family of tori of revolution. 
Since the energy we use is conformally invariant, it follows that the Clifford torus gives the minimum energy among Dupin cyclides. 
As the surface energy that we use in this paper is generelization of knot energy, we start with the review of it. 

{\em Energy of knots} was introduced in \cite{O} motivated to give a functional on the space of knots that can produce a representative configuration of a knot for each knot type as an embedding that minimizes the energy in the knot type. 
Let $K$ be a knot and $x$ be a point on it. 
Define 
\begin{equation}\label{def_V_E}
\begin{array}{rcl}
V(x;K)&=&\displaystyle \lim_{\varepsilon\to0}\left(\int_{K\setminus B_\varepsilon(x)}\frac{dy}{|x-y|^2}-\frac2\varepsilon\right), \>\>\mbox{ and }\>\> 
E(K)
=\displaystyle \int_KV(x;K)dx,
\end{array}
\end{equation}
where  $B_\varepsilon(x)$ is a ball with center $x$ and radius $\varepsilon$. 
%
Let us call a process as in the definition of $V(x;K)$ in \eqref{def_V_E} the {\em renormalization} in this paper. 
In general, when we are interested in a diverging integral, we first restrict the integration to the complement of an $\varepsilon$-neighbourhood of the set where the integrand blows up, then expand the result in a Laurent series in $\varepsilon$, and finally take the constant term. 
In the case of a knot, the integrand of $V$ in \eqref{def_V_E} blows up at the one-point set $\{x\}$. There are two ways to define an $\varepsilon$-neighbourhood of it, according to the choice of the distance between a pair of points on the knot; either the chord length as in \eqref{def_V_E} or the arc-length along the knot as in \cite{O}. Both types of the renormalization give the same result (\cite{Obook}). 

The energy $E(K)$ in \eqref{def_V_E} was proved to be invariant under M\"obius transformations by Freedman, He, and Wang (\cite{FHW}), which is the reason why it is sometimes called the {\em M\"obius energy} of knots. 

After this energy was found, it has been generalized to functionals that can measure geometric complexity of knots, surfaces, and in general, submanifolds (\cite{AS}, \cite{KS}, et al.). 
Among several ways of generalization to surface energy, in this paper we study the one by Auckly and Sadun that uses a similar renormalization process as in \eqref{def_V_E}. 

Let $S$ be an embedded surface in $\mathbb{R}^3$ without boundary and $x$ be a point in $S$. 
Define the {\em renormalized $r^{-4}$-potential} $V$ and the {\em renormalized $r^{-4}$-potential energy} $E$ by 
\begin{equation}\label{def_V_E_surface}
\begin{array}{rcl}
V(x;S)&=&\displaystyle \lim_{\varepsilon \to 0} \left(\int_{S\setminus B_\varepsilon(x)}\frac{d^2y}{|x-y|^4}-\frac{\pi}{\varepsilon^2}+\frac{\pi \Delta(x)}{16}\log\left(\Delta(x)\varepsilon^2\right)+\frac{\pi K(x)}{4}\right), \\[4mm]
E(S)&=&\displaystyle \int_SV(x;S)d^2x,
\end{array}
\end{equation}
where $d^2y$ and $d^2x$ mean the volume element of $S$, $\Delta(x)$ is given by 
$\Delta(x)
=(\kappa_1(x)-\kappa_2(x))^2,  
$
%
and $K(x)$ is the Gauss curvature; $K=\kappa_1\kappa_2$. 
%
This energy $E(S)$ was proved to be invariant under M\"obius transformations in \cite{AS}. 
It blows up as $S$ degenerates to an immersed surface with double points. 

As was pointed out in \cite{AS}, the choice of the $\log$-term in \eqref{def_V_E_surface} is not the unique reasonable one. 
The reason why there is a factor $\Delta(x)$ in the $\log$-term is to make the resulting energy scale invariant, but $c\Delta(x)$ $(c\ne0)$ also has the same effect. Thus there is ambiguity in the definition of the renormalized potential. 

\smallskip
In \cite{AS}, Auckly and Sadun has computed the energy of spheres and planes and the potentials $V$ of an infinitely long straight cylinder and a surface called dimple. 
In this article we compute the energies of one-parameter family of tori of revolution. 

\bigskip\noindent
{\bf Theorem}: {\em 
Let $T_R$ be a torus of revolution whose generating circle has radius $1$ and center at distance $R$ $(R>1)$ from the axis of revolution. 
Then the renormalized potential energy is given by 
\[
E(T_{R})=\frac{\pi^3}{2\sqrt{R^2-1}} \left({R^2} \left( 3\log 2 -1 \right) +2-\frac{2}{R^2}\right). 
\]
}


\bigskip
\noindent
{\bf Corollary}: {\em Among tori of revolution, a stereographic projection of the Clifford torus gives the minimum energy. 
}

\bigskip\noindent
{\bf Proof of Corollary}: Since
\[
\frac{d}{dR}E(T_{R})
=\frac {\pi^3 ({R}^{2}-2)  \left(({R}^{2}-2) ^{2}+3\,R^4\log  \left( 4/e\right)  \right) }
{4 R^3\left(R^2-1 \right) ^{3/2}},
\]
the energy takes the minimum value $\pi^3(6\log2-1)/2$ when $R=\sqrt2$. {\hfill$\square$\par\medskip}

\bigskip
\noindent
{\bf Problem}:
(1) Does $T_{\sqrt2}$ give the minimum energy among all the embedded tori in $\mathbb R^3$, hence among all the immersed tori in $\mathbb R^3$?

(2) When we change the power of $|x-y|$ in the denominator in \eqref{def_V_E_surface} from $4$ to any number $\lambda$, we obtain a new potential energy $E_{r^{-\lambda}}$ after suitable renormalization. It is no longer scale invariant when $\lambda\ne4$. What is $R=R(\lambda)$ that makes $T_R$ give the minimum energy $E_{r^{-\lambda}}$ after rescaling to have area $1$?

\section{Computation of the energy of a torus}

\bigskip\noindent
{\bf Proof of Theorem}: Let $T$ be a torus of revolution parametrized by 
\[
p(u,v)=((R+\cos u)\cos v,(R+\cos u)\sin v,\sin u) 
.
\]

We compute the renormalized $r^{-4}$-potential of $T$ at a point $x=p(\alpha,0)=(R+\cos \alpha,0,\sin \alpha)$. 

\medskip
Let $\textrm{Dist}=\textrm{Dist}(u,v)$ be the distance between $x$ and a point $y=p(u,v)$:
\begin{eqnarray*}
\textrm{Dist}^2&=&|x-y|^2 
=|p(u,v)-p(\alpha,0)|^2\\
&=& \left( \left( R+\cos\alpha\right)-\left( R+\cos u \right)\cos v\right)^2 +\left( R+\cos u \right)^2\sin^2 v+\left( \sin\alpha -\sin u \right)^2\\[1mm]
   &=& 2R^2+2+2R\left( \cos\alpha +\cos u \right)-2\sin\alpha\sin u-2\left(R+\cos\alpha\right)\left(R+\cos u\right) \cos v.
\end{eqnarray*}
First put 
$$t=2\sin \frac{u-\alpha}{2},\ \ s=2\left( R+\cos\alpha\right) \sin \frac{v}{2},$$
then, as we have 
\[\begin{array}{c}
\cos u =\displaystyle \frac{2-t^2}{2}\cos\alpha-\frac{t\sqrt{4-t^2}}{2}\sin\alpha, \ \ 
\sin u =\displaystyle \frac{t\sqrt{4-t^2}}{2}\cos\alpha+\frac{2-t^2}{2}\sin\alpha, \\[2mm]
\cos v = \displaystyle \frac{2\left(R+\cos\alpha\right)^2-s^2}{2\left(R+\cos\alpha\right)^2},
\end{array}\]
the distance can be expressed as 
\begin{eqnarray*}
\textrm{Dist}^2 &=& t^2+s^2-\frac{\cos\alpha}{2\left (R+\cos\alpha\right)}t^2s^2-\frac{\sin\alpha}{2\left(R+\cos\alpha\right)}s^2t\sqrt{4-t^2}. 
\end{eqnarray*}
Next, put 
\[
\theta=\frac{\pi+\alpha-u}2, \ \ \varphi=\frac{\pi-v}2, 
\] 
then, as 
\[t=2\cos \theta,\ \ s=2(R+\cos\alpha)\cos \varphi,\]
the distance can be expressed as 
\[
\textrm{Dist}^2=\cos^2\theta+
(R+\cos\alpha)\left(R+\cos\alpha-2\cos\alpha\cos^2 \theta-2\sin\alpha|\sin\theta|\cos\theta\right)\cos^2\varphi.
\]

\medskip
Put 
\[
V(\varepsilon,x)=\iint_{\textrm{Dist}\ge\varepsilon} \frac{d^2y}{|x-y|^2}.
\]
Since the area element of $T$ is given by $d^2y=(R+\cos u)\,dudv$, $V(\varepsilon,x)$ is given by 
\begin{eqnarray*}
&& 
%
\iint_{\textrm{Dist}\ge\varepsilon} \frac
{(R+\cos u)}
{\left(2R^2+2+2R\left( \cos\alpha +\cos u \right)-2\sin\alpha\sin u-2\left(R+\cos\alpha\right)\left(R+\cos u\right) \cos v\right)^2}
\,dudv \\[2mm]
  &=& \iint_{\textrm{Dist}\ge\varepsilon} 
\frac{ \left(R+\left(\frac{2-t^2}{2}\cos\alpha-\frac{t\sqrt{4-t^2}}{2}\sin\alpha\right)\right)}
{ \left(t^2+s^2-\frac{\cos\alpha}{2\left (R+\cos\alpha\right)}t^2s^2-\frac{\sin\alpha}{2\left(R+\cos\alpha\right)}s^2t\sqrt{4-t^2}\right)^2} \\
&&\hspace{5cm}
\times\frac{2}{\sqrt{4-t^2}}  \frac{2}{\sqrt{4(R+\cos\alpha)^2-s^2}} \, dtds \\[2mm]
&=& \frac{r}{4}\iint_{\textrm{Dist}\ge\varepsilon} \frac
{\left(R+\cos\alpha-2\cos\alpha\cos^2 \theta-2\sin\alpha|\sin\theta|\cos\theta\right)\cdot\textrm{sgn}(\sin\theta)\cdot\textrm{sgn}(\sin\varphi)}
{\left(\cos^2\theta+
(R+\cos\alpha)\left(R+\cos\alpha-2\cos\alpha\cos^2 \theta-2\sin\alpha|\sin\theta|\cos\theta\right)\cos^2\varphi\right)^2}
\, d\theta d\varphi.
%
\end{eqnarray*}

\medskip
We devide the half of the above integral into four parts, $I_1,I_2,I_3$, and $I_4$, which are the integrals over the following four regions respectively:
\[\begin{array}{l}
\left\{p(u,v)\,|\,\alpha+2\sin^{-1}(\varepsilon/2 )\le u\le \alpha+\pi, 0\le v\le \pi\right\},\\[2mm]
\left\{p(u,v)\,|\,\alpha \le u \le \alpha+2\sin^{-1}(\varepsilon/2 ), 0\le v\le \pi, |p(u,v)-x|\ge\varepsilon\right\},\\[2mm]
\left\{p(u,v)\,|\,\alpha-\pi\le u \le \alpha-2\sin^{-1}(\varepsilon/2 ), 0\le v\le \pi\right\},\\[2mm]
\left\{p(u,v)\,|\,\alpha-2\sin^{-1}(\varepsilon/2 )\le u\le \alpha, 0\le v\le \pi, |p(u,v)-x|\ge\varepsilon\right\}.
\end{array}\]

\medskip
Let us first compute $I_2$ and $I_4$. 
Put 
\[
c(t)= 4(R+\cos\alpha)^2-{2(R+\cos\alpha)\cos\alpha}\cdot t^2-{2(R+\cos\alpha)\sin\alpha}{r}\cdot t\sqrt{4-t^2}.
\]
Define $\vartheta$ by $t=\varepsilon\cos\vartheta$. 
Then we have 
\begin{eqnarray*}
\begin{array}{rcl}
I_2 &=&\displaystyle  \int^{\varepsilon}_{0} \int^{\pi}_{2\sin^{-1}\left(\sqrt{\frac{\varepsilon^2-t^2}{c(t)}}\right)} \frac{2R+\cos\alpha\cdot(2-t^2)-\sin\alpha\cdot t\sqrt{4-t^2}}{(t^2+c(t)\sin^2\frac{v}{2})^2}\frac{dt}{\sqrt{4-t^2}}\,dv \\[5mm]
  &=& \displaystyle \int^{\varepsilon}_{0} \frac{2R+\cos\alpha\cdot(2-t^2)-\sin\alpha\cdot t\sqrt{4-t^2}}{\sqrt{4-t^2}} \\ [4mm]
  & & \displaystyle \times  \left[ \frac{2t^2+c(t)}{t^3(t^2+c(t))^\frac{3}{2}} \tan^{-1}\left( \frac{\sqrt{t^2+c(t)}}{t} \tan\frac{v}{2}\right) \right. \\[4mm]
  & & \displaystyle \hspace{2cm} \left.
+\frac{c(t)\tan\frac{v}{2}}{t^2(t^2+c(t))\left((t^2+c(t))\tan^2 \frac{v}{2}+t^2\right)}\right]^{\pi}_{v=2\sin^{-1}\left(\sqrt{\frac{\varepsilon^2-t^2}{c(t)}}\right)} dt \\
  &=& \displaystyle  \int^{\varepsilon}_{0} \frac{2R+\cos\alpha\cdot(2-t^2)-\sin\alpha\cdot t\sqrt{4-t^2}}{\sqrt{4-t^2}} \\[4mm]
  & &\displaystyle \times  \left\{ \frac{2t^2+c(t)}{t^3(t^2+c(t))^\frac{3}{2}} \left( \frac{\pi}{2}-\tan^{-1}\left( \frac{\sqrt{\varepsilon^2-t^2}}{t}\sqrt{\frac{c(t)+t^2}{c(t)+t^2-\varepsilon^2}} \right)\right) -\frac{\sqrt{\varepsilon^2-t^2}\sqrt{c(t)+t^2-\varepsilon^2}}{\varepsilon^2t^2(t^2+c(t))}\right\}dt\\[5mm]
  &=& \displaystyle \int^{\frac{\pi}{2}}_{0} \frac{2R+\cos\alpha\cdot(2-\varepsilon^2\cos^2\vartheta)-\varepsilon\sin\alpha\cos\vartheta\sqrt{4-\varepsilon^2\cos^2\vartheta}}{\sqrt{4-\varepsilon^2\cos^2\vartheta}} \\[4mm]
  & & \hspace{0.2cm} \displaystyle \times \left\{ \frac{\sin\vartheta\left(2\varepsilon^2\cos^2\vartheta+c(\varepsilon\cos\vartheta)\right)}{\varepsilon^2\cos^3\vartheta(\varepsilon^2\cos^2\vartheta+c(\varepsilon\cos\vartheta))^\frac{3}{2}} \left( \frac{\pi}{2}-\tan^{-1}\left( \tan\vartheta\sqrt{\frac{c(\varepsilon\cos\vartheta)+\varepsilon^2\cos^2\vartheta}{c(\varepsilon\cos\vartheta)+\varepsilon^2\cos^2\vartheta-\varepsilon^2}} \right)\right) \right. \\[6mm]
 & & \hspace{0.8cm} \left. 
-\displaystyle \frac{\tan^2\vartheta \sqrt{c(\varepsilon\cos\vartheta)+\varepsilon^2\cos^2\vartheta-\varepsilon^2}}{\varepsilon^2(\varepsilon^2\cos^2\vartheta+c(\varepsilon\cos\vartheta))}\,\right\}\,d\vartheta. \\[5mm]
\end{array}
\end{eqnarray*}
Expand the integrand in a series in $\varepsilon$ to obtain 
\begin{eqnarray*}
\begin{array}{rcl}
I_2 &=& \displaystyle \frac{1}{2\varepsilon^2}\int^{\frac{\pi}{2}}_{0}\left( \left(\frac{\pi}{2}-\vartheta \right)\frac{\sin\vartheta}{\cos^3\vartheta}-\frac{\sin^2\vartheta}{\cos^2\vartheta} \right)d\vartheta
-\frac{\sin\alpha}{4(R+\cos\alpha)\varepsilon}\int^{\frac{\pi}{2}}_{0}\left( \left(\frac{\pi}{2}-\vartheta \right)\frac{\sin\vartheta}{\cos^2\vartheta}-\frac{\sin^2\vartheta}{\cos\vartheta} \right)d\vartheta \\[4mm]
  & &\displaystyle  +\frac{R^2}{16(R+\cos\alpha)^2}\int^{\frac{\pi}{2}}_{0}\left( \left(\frac{\pi}{2}-\vartheta \right)\frac{\sin\vartheta}{\cos\vartheta}-\sin^2\vartheta \right)d\vartheta  
+\frac{1}{8(R+\cos\alpha)^2}\int^{\frac{\pi}{2}}_{0} \sin^2\vartheta d\vartheta +O(\varepsilon).
\end{array}
\end{eqnarray*}
Since 
\[
\begin{array}{c}
\displaystyle \int^{\frac{\pi}{2}}_{0}\left( \left(\frac{\pi}{2}-\vartheta \right)\frac{\sin\vartheta}{\cos^3\vartheta}-\frac{\sin^2\vartheta}{\cos^2\vartheta} \right)d\vartheta=\frac\pi4, \ \ 
\int^{\frac{\pi}{2}}_{0}\left( \left(\frac{\pi}{2}-\vartheta \right)\frac{\sin\vartheta}{\cos^2\vartheta}-\frac{\sin^2\vartheta}{\cos\vartheta} \right)d\vartheta=\frac{4-\pi}2,\\[4mm]
\displaystyle \int^{\frac{\pi}{2}}_{0}\left( \left(\frac{\pi}{2}-\vartheta \right)\frac{\sin\vartheta}{\cos\vartheta}-\sin^2\vartheta \right)d\vartheta =\frac{\pi}4(2\log2-1),
\end{array}
\]
we have
\begin{eqnarray*}
\begin{array}{rcl}
I_2 &=&\displaystyle \frac{\pi}{8\varepsilon^2}+\frac{(\pi-4)\sin\alpha}{8(R+\cos\alpha)\varepsilon}+\frac{\pi R^2}{64(R+\cos\alpha)^2}(2\log 2-1)+\frac{\pi}{32(R+\cos\alpha)^2}+O(\varepsilon).
\end{array}
\end{eqnarray*}
The integral $I_4$ can be obtained from $I_2$ by changing $\alpha$ to $-\alpha$: 
\begin{eqnarray*}
I_4 &=& \int^{0}_{-\varepsilon} \int^{\pi}_{2\sin^{-1}\left(\sqrt{\frac{\varepsilon^2-t^2}{c(t)}}\right)} \frac{2R+\cos\alpha\cdot(2-t^2)-\sin\alpha\cdot t\sqrt{4-t^2}}{(t^2+c(t)\sin^2\frac{v}{2})^2}\frac{dt}{\sqrt{4-t^2}}\,dv \\
   &=& \frac{\pi}{8\varepsilon^2}-\frac{(\pi-4)\sin\alpha}{8(R+\cos\alpha)\varepsilon}+\frac{\pi R^2}{64(R+\cos\alpha)^2}(2\log 2-1)+\frac{\pi}{32(R+\cos\alpha)^2}+O(\varepsilon).
\end{eqnarray*}

\medskip
Let us next compute $I_1$ and $I_3$. 
The integral $I_1$ is given by 
\[
I_1=\frac{r}{4} \int^{\cos^{-1}\frac{\varepsilon}{2}}_{0} \!\!\!\int^{\frac{\pi}{2}}_{0}\frac{R+\cos\alpha-2\cos\alpha\cos^2 \theta-2\sin\alpha\sin\theta\cos\theta}{(\cos^2 \theta +(R+\cos\alpha)(R+\cos\alpha-2\cos\alpha \cos^2 \theta -2\sin\alpha\sin\theta\cos\theta)\cos^2 \varphi)^2}\, d\theta d\varphi.
\]
Put 
\[a=\cos\theta,\ \  b=(R+\cos\alpha)\left(R+\cos\alpha-2\cos\alpha\cos^2 \theta-2\sin\alpha\sin\theta\cos\theta\right)\]
so that the denominator of the integrand is given by 
$\textrm{Dist}^4=(a^2+b\,\cos^2\varphi)^2$.
Then 
\begin{eqnarray*}
I_1 &=& \displaystyle \frac{r}{4} \int^{\cos^{-1}\frac{\varepsilon}{2}}_{0} \left(R+\cos\alpha-2\cos\alpha\cos^2 \theta-2\sin\alpha\sin\theta\cos\theta\right) \\
   & & \displaystyle \times \left[\frac{2a^2+b}{2a^3(a^2+b)^\frac{3}{2}}\tan^{-1}\left( \frac{a}{\sqrt{a^2+b}}\tan\varphi \right) -\frac{b\tan\varphi}{2a^2(a^2+b)(a^2\tan^2\varphi+a^2+b)}\right]^{\frac{\pi}{2}}_{\varphi=0}d\theta\\
   &=& \frac{r}{4} \int^{\cos^{-1}\frac{\varepsilon}{2}}_{0} \left(R+\cos\alpha-2\cos\alpha\cos^2 \theta-2\sin\alpha\sin\theta\cos\theta\right) \frac{(a^2+b)+a^2}{2a^3(a^2+b)^\frac{3}{2}}\cdot \frac{\pi}{2} \,d\theta \\
   &=&\frac{\pi(R+\cos\alpha)}{16}I_{11}-\frac{\pi\cos\alpha}{8}I_{12}-\frac{\pi\sin\alpha}{8}I_{13}+\frac{\pi(R+\cos\alpha)}{16}I_{14}-\frac{\pi \cos\alpha}{8}I_{15}-\frac{\pi \sin\alpha}{8}I_{16}, 
\end{eqnarray*}
where
\[
\begin{array}{c}
\displaystyle I_{11}=\int^{\cos^{-1}\frac{\varepsilon}{2}}_{0}\frac{\,d\theta}{a^3\sqrt{a^2+b}}, \ \
I_{12}=\int^{\cos^{-1}\frac{\varepsilon}{2}}_{0}\frac{\cos^2\theta\,d\theta}{a^3\sqrt{a^2+b}}, \ \
I_{13}=\int^{\cos^{-1}\frac{\varepsilon}{2}}_{0}\frac{\sin\theta\cos\theta\,d\theta}{a^3\sqrt{a^2+b}},\\[4mm]
\displaystyle I_{14}=\int^{\cos^{-1}\frac{\varepsilon}{2}}_{0}\frac{r\,d\theta}{a(a^2+b)^{\frac32}}, \ \
I_{15}=\int^{\cos^{-1}\frac{\varepsilon}{2}}_{0}\frac{\cos^2\theta\,d\theta}{a(a^2+b)^{\frac32}}, \ \
I_{16}=\int^{\cos^{-1}\frac{\varepsilon}{2}}_{0}\frac{\sin\theta\cos\theta\,d\theta}{a(a^2+b)^{\frac32}}.
\end{array}
\]
Since 
\[
a^2+b=\cos^2\theta\left(R^2+((R+\cos\alpha)\tan\theta-\sin\alpha)^2\right), \ \ 
\tan\left(\cos^{-1}\left(\frac\varepsilon{2}\right)\right)=\frac{\sqrt{4-\varepsilon^2}}\varepsilon,
\]
we have
\begin{eqnarray*}
I_{15} &=& \int^{\cos^{-1}\frac{\varepsilon}{2}}_{0} \frac{1}{\cos^2\theta\left(R^2+\left((R+\cos\alpha)\tan\theta-\sin\alpha\right)^2\right)^{\frac{3}{2}}}\, d\theta\\
  &=& \left[ \frac{(R+\cos\alpha)\tan\theta-\sin\alpha}{R^2(R+\cos\alpha)\sqrt{R^2+\left((R+\cos\alpha)\tan\theta-\sin\alpha\right)^2}}\right]^{\cos^{-1}\frac{\varepsilon}{2}}_{0}\\
  &=&\frac{(R+\cos\alpha)\sqrt{4-\varepsilon^2}-(\sin\alpha)\varepsilon}{R^2(R+\cos\alpha)\sqrt{R^2\varepsilon^2+((R+\cos\alpha)\sqrt{4-\varepsilon^2}-(\sin\alpha)\varepsilon)^2}}
+\frac{\sin\alpha}{R^2(R+\cos\alpha)\sqrt{R^2+\sin^2\alpha}}\\
  &=& \frac{1}{R^2(R+\cos\alpha)}+\frac{\sin\alpha}{R^2(R+\cos\alpha)\sqrt{R^2+\sin^2\alpha}}+O\left(\varepsilon\right),
\end{eqnarray*}

\begin{eqnarray*}
I_{16} &=& \int^{\cos^{-1}\frac{\varepsilon}{2}}_{0} \frac{\sin\theta}{\cos^3\theta\left(R^2+\left((R+\cos\alpha)\tan\theta-\sin\alpha\right)^2\right)^{\frac{3}{2}}}\, d\theta\\
  &=& \left[ \frac{\sin\alpha\left((R+\cos\alpha)\tan\theta-\sin\alpha\right)-R^2}{R^2(R+\cos\alpha)^2\sqrt{R^2+\left((R+\cos\alpha)\tan\theta-\sin\alpha\right)^2}}\right]^{\cos^{-1}\frac{\varepsilon}{2}}_{0}\\
  &=&\frac{\sin\alpha\left((R+\cos\alpha)\sqrt{4-\varepsilon^2}-(\sin\alpha)\varepsilon\right)-R^2\varepsilon}{R^2(R+\cos\alpha)^2\sqrt{R^2\varepsilon^2+\left((R+\cos\alpha)\sqrt{4-\varepsilon^2}-(\sin\alpha)\varepsilon\right)^2}}
+\frac{\sqrt{R^2+\sin^2\alpha}}{R^2(R+\cos\alpha)^2} \\
  &=& \frac{\sin\alpha}{R^2(R+\cos\alpha)^2}+ \frac{\sqrt{R^2+\sin^2\alpha}}{R^2(R+\cos\alpha)^2}+O\left(\varepsilon\right),
\end{eqnarray*}

\begin{eqnarray*}
I_{12} &=& \int^{\cos^{-1}\frac{\varepsilon}{2}}_{0} \frac{1}{\cos^2\theta\sqrt{R^2+\left((R+\cos\alpha)\tan\theta-\sin\alpha\right)^2}}\, d\theta\\
  &=& \left[ \frac{1}{R+\cos\alpha}\sinh^{-1}\left(\frac{(R+\cos\alpha)\tan\theta-\sin\alpha}{R}\right)\right]^{\cos^{-1}\frac{\varepsilon}{2}}_{0}\\
  &=& \frac{1}{R+\cos\alpha} \left\{  \sinh^{-1}\left(\frac{(R+\cos\alpha)\sqrt{4-\varepsilon^2}-(\sin\alpha)\varepsilon}{R\varepsilon}\right)+\sinh^{-1}\left(\frac{\sin\alpha}{R}\right)\right\} \\
  &=& \frac{1}{R+\cos\alpha} \left\{  \log\left(\frac{4(R+\cos\alpha)}{R\varepsilon}\right)+O\left(\varepsilon\right)+\sinh^{-1}\left(\frac{\sin\alpha}{R}\right)\right\}\\
  &=& -\frac{1}{2(R+\cos\alpha)}\log\left( \frac {R^2\varepsilon^2}{16(R+\cos\alpha)^2} \right) +\frac{1}{R+\cos\alpha}\sinh^{-1}\left(\frac{\sin\alpha}{R}\right)+O(\varepsilon),
\end{eqnarray*}

\begin{eqnarray*}
I_{13} &=& \int^{\cos^{-1}\frac{\varepsilon}{2}}_{0} \frac{\sin\theta}{\cos^3\theta\sqrt{R^2+\left((R+\cos\alpha)\tan\theta-\sin\alpha\right)^2}}\, d\theta \\
 &=& \frac{\sin\alpha}{R+\cos\alpha}I_{12}+\frac{1}{(R+\cos\alpha)^2}\left[\sqrt{R^2+\left((R+\cos\alpha)\tan\theta-\sin\alpha\right)^2}\> \right]^{\cos^{-1}\frac{\varepsilon}{2}}_{0}  \\[1mm]
  &=& \frac{\sin\alpha}{R+\cos\alpha}I_{12} +\frac{\sqrt{R^2\varepsilon^2+\left((R+\cos\alpha)\sqrt{4-\varepsilon^2}-(\sin\alpha)\varepsilon\right)^2}}{(R+\cos\alpha)^2\varepsilon}
-\frac{\sqrt{R^2+\sin^2\alpha}}{(R+\cos\alpha)^2} \\[1mm]
 &=& \frac{\sin\alpha}{R+\cos\alpha}I_{12} +\left\{  \frac{2}{(R+\cos\alpha)\varepsilon}-\frac{\sin\alpha}{(R+\cos\alpha)^2}+O\left( \varepsilon\right) \right\}-\frac{\sqrt{R^2+\sin^2\alpha}}{(R+\cos\alpha)^2}\,,
\end{eqnarray*}

\begin{eqnarray*}
I_{14} &=& \int^{\cos^{-1}\frac{\varepsilon}{2}}_{0} \frac{1}{\cos^4\theta\left(R^2+\left((R+\cos\alpha)\tan\theta-\sin\alpha\right)^2\right)^{\frac{3}{2}}}\, d\theta\\
  &=& \frac{1}{(R+\cos\alpha)^2}I_{12}-\frac{(1-2R\cos\alpha-2\cos^2\alpha)}{(R+\cos\alpha)^2}I_{15}+\frac{2\sin\alpha}{R+\cos\alpha}I_{16} \\
  &=& \frac{1}{(R+\cos\alpha)^2}I_{12} +\frac{1+2R\cos\alpha}{R^2(R+\cos\alpha)^3} +\frac{(1+2R^2)\sin\alpha+2R^2\sin\alpha\cos\alpha}{R^2(R+\cos\alpha)^3 \sqrt{R^2+\sin^2\alpha}}+O(\varepsilon)\,,
\end{eqnarray*}


and
\begin{eqnarray*}
I_{11} &=& \int^{\cos^{-1}\frac{\varepsilon}{2}}_{0} \frac{1}{\cos^4\theta\sqrt{R^2+\left((R+\cos\alpha)\tan\theta-\sin\alpha\right)^2}}\, d\theta\\
  &=&  \frac{R^2-1+4R\cos\alpha+3\cos^2\alpha}{2(R+\cos\alpha)^2}I_{12}+\frac{3\sin\alpha}{2(R+\cos\alpha)}I_{13}\\
  && +\left[ \frac{\sin\theta\sqrt{R^2+\left((R+\cos\alpha)\tan\theta-\sin\alpha\right)^2}}{2(R+\cos\alpha)^2\cos\theta}\right]^{\cos^{-1}\frac{\varepsilon}{2}}_{0}\\
  &=& \frac{R^2-1+4R\cos\alpha+3\cos^2\alpha}{2(R+\cos\alpha)^2}I_{12}+\frac{3\sin\alpha}{2(R+\cos\alpha)}I_{13} \\
  && + \frac{\sqrt{4-\varepsilon^2}\sqrt{R^2\varepsilon^2+\left((R+\cos\alpha)\sqrt{4-\varepsilon^2}-(\sin\alpha)\varepsilon\right)^2}}{2(R+\cos\alpha)^2\varepsilon^2}\\
  &=&   \frac{R^2-1+4R\cos\alpha+3\cos^2\alpha}{2(R+\cos\alpha)^2}I_{12}+\frac{3\sin\alpha}{2(R+\cos\alpha)}I_{13} \\[1mm]
  && + \left\{  \frac {2}{(R+\cos\alpha)\varepsilon^2} -\frac{\sin\alpha}{(R+\cos\alpha)^2\varepsilon}-\frac{R^2+4R\cos\alpha+2\cos^2\alpha}{4(R+\cos\alpha)^3} +O\left(\varepsilon \right)  \right\} \,.
\end{eqnarray*}
%
Therefore
\begin{eqnarray*}
\begin{array}{lll}
I_1&=&\displaystyle \frac{\pi(R+\cos\alpha)}{16}I_{11}-\frac{\pi\cos\alpha}{8}I_{12}-\frac{\pi\sin\alpha}{8}I_{13}+\frac{\pi(R+\cos\alpha)}{16}I_{14}-\frac{\pi \cos\alpha}{8}I_{15}-\frac{\pi \sin\alpha}{8}I_{16} \\[4mm]
   &=& \displaystyle \frac{\pi R^2}{32(R+\cos\alpha)}I_{12}+\frac{\pi}{8\varepsilon^2}-\frac{\sin\alpha}{8(R+\cos\alpha)\varepsilon}-\frac{\pi(R^2+4R\cos\alpha+2\cos^2\alpha)}{64(R+\cos\alpha)^2}-\frac{\pi }{16R^2(R+\cos\alpha)^2}\\ [4mm]
  & & \displaystyle +\frac{\pi\sin^2\alpha}{32(R+\cos\alpha)^2}-\frac{\pi \sin\alpha}{16R^2(R+\cos\alpha)^2\sqrt{R^2+\sin^2\alpha}}+\frac{\pi\sin\alpha \sqrt{R^2+\sin^2\alpha}}{32(R+\cos\alpha)^2}+O(\varepsilon)\\[4mm]
  &=& \displaystyle -\frac{\pi R^2}{64(R+\cos\alpha)^2}\log\left( \frac {R^2\varepsilon^2}{16(R+\cos\alpha)^2} \right) + \frac{\pi R^2}{32(R+\cos\alpha)^2} \sinh^{-1}\left( \frac{\sin\alpha}{R}\right)\\ [4mm]
& & \displaystyle  +\frac{\pi}{8\varepsilon^2}-\frac{\sin\alpha}{8(R+\cos\alpha)\varepsilon}-\frac{\pi(R^2+4R\cos\alpha+2\cos^2\alpha)}{64(R+\cos\alpha)^2}-\frac{\pi }{16R^2(R+\cos\alpha)^2}+\frac{\pi\sin^2\alpha}{32(R+\cos\alpha)^2}
\\ [4mm]
& & \displaystyle  -\frac{\pi \sin\alpha}{16R^2(R+\cos\alpha)^2\sqrt{R^2+\sin^2\alpha}}+\frac{\pi\sin\alpha \sqrt{R^2+\sin^2\alpha}}{32(R+\cos\alpha)^2}+O(\varepsilon).
\end{array}
\end{eqnarray*}
The integral $I_3$ can be obtained from $I_1$ by changing $\alpha$ to $-\alpha$:
\begin{eqnarray*}
I_3 &=& \frac{r}{4} \int^{\cos^{-1}\frac{\varepsilon}{2}}_{0} \!\!\!\int^{\frac{\pi}{2}}_{0}\frac{R+\cos\alpha-2\cos\alpha\cos^2 \theta+2\sin\alpha\sin\theta\cos\theta}{(\cos^2 \theta +(R+\cos\alpha)(R+\cos\alpha-2\cos\alpha \cos^2 \theta +2\sin\alpha\sin\theta\cos\theta)\cos^2 \varphi)^2}\, d\theta d\varphi  \\
  &=& \displaystyle -\frac{\pi R^2}{64(R+\cos\alpha)^2}\log\left( \frac {R^2\varepsilon^2}{16(R+\cos\alpha)^2} \right) - \frac{\pi R^2}{32(R+\cos\alpha)^2} \sinh^{-1}\left( \frac{\sin\alpha}{R}\right)\\ 
& & \displaystyle  +\frac{\pi}{8\varepsilon^2}+\frac{\sin\alpha}{8(R+\cos\alpha)\varepsilon}-\frac{\pi(R^2+4R\cos\alpha+2\cos^2\alpha)}{64(R+\cos\alpha)^2}-\frac{\pi }{16R^2(R+\cos\alpha)^2}+\frac{\pi\sin^2\alpha}{32(R+\cos\alpha)^2}
\\ 
& & \displaystyle  +\frac{\pi \sin\alpha}{16R^2(R+\cos\alpha)^2\sqrt{R^2+\sin^2\alpha}}=\frac{\pi\sin\alpha \sqrt{R^2+\sin^2\alpha}}{32(R+\cos\alpha)^2}+O(\varepsilon).
%
\end{eqnarray*}

\medskip
By putting all the formulae together we obtain 
\begin{eqnarray*}
V\left(\varepsilon,x\right)&=& 2\left(I_1+I_2+I_3+I_4\right)\\
  &=&  \frac{\pi}{\varepsilon^2}-\frac{\pi R^2}{16(R+\cos\alpha)^2}\log\left( \frac {R^2\varepsilon^2}{(R+\cos\alpha)^2} \right) +\frac{\pi R^2}{8(R+\cos\alpha)^2}3\log 2\\
  &&-\frac{\pi}{8}-\frac{\pi }{4R^2(R+\cos\alpha)^2}+\frac{\pi(1+\sin^2\alpha)}{8(R+\cos\alpha)^2}.
\end{eqnarray*}

As the Gauss curvature and $\Delta=(\kappa_1-\kappa_2)^2$ at the point $x=p(\alpha,0)$ is given by
\[K(x)=\frac{\cos \alpha}{\left(R+\cos \alpha\right)}, \ \ \Delta(x)=\frac{R^2}{\left( R+\cos \alpha\right)^2}\,,\]
the renormalized potential is given by 
\begin{eqnarray*}
V(x;T)&=&\lim_{\varepsilon \to 0} \left(V\left(\varepsilon,x\right)-\frac{\pi}{\varepsilon^2}+\frac{\pi \Delta(x)}{16}\log\left(\Delta(x)\varepsilon^2\right)+\frac{\pi K(x)}{4}\right) \\
&=& \frac{\pi R^2}{8(R+\cos\alpha)^2}3\log 2-\frac{\pi}{8}-\frac{\pi }{4R^2(R+\cos\alpha)^2}+\frac{\pi(1+\sin^2\alpha)}{8(R+\cos\alpha)^2}+\frac{\pi\cos\alpha}{4(R+\cos\alpha)}.
\end{eqnarray*}
It implies that the renormalized $r^{-4}$-potential energy of the torus $T$ is given by 
\begin{eqnarray*}
E (T)&=& \int_T V(x;T)\, d^2x \\
  &=&  2\pi \int^{2\pi}_{0} \left(\frac{\pi R^2}{8(R+\cos\alpha)^2}3\log 2-\frac{\pi}{8}-\frac{\pi }{4R^2(R+\cos\alpha)^2}+\frac{\pi(1+\sin^2\alpha)}{8(R+\cos\alpha)^2}\right. \\
  && \left. \hspace{1.5cm}
+\frac{\pi\cos\alpha}{4(R+\cos\alpha)} \right)(R+\cos\alpha)\,d\alpha\\
  &=& \frac{\pi^2 r}{2} \int^{\pi}_{0} \left(\left({R^2} \left( 3\log 2 -1 \right) +2-\frac{2}{R^2}\right) \frac{1}{R+\cos\alpha} \right) d\alpha \\
  &=& \frac{\pi^2 r}{2} \left[ \left({R^2} \left( 3\log 2 -1 \right) +2-\frac{2}{R^2}\right)\frac{2}{\sqrt{R^2-1}}\tan^{-1}\left(\sqrt{\frac{R-1}{R+1}}\tan\frac{\alpha}{2}\right) \right]^{\pi}_{0} \\
  &=&\frac{\pi^3r}{2\sqrt{R^2-1}} \left({R^2} \left( 3\log 2 -1 \right) +2-\frac{2}{R^2}\right).
\end{eqnarray*}

{\hfill$\square$\par\medskip} 
\section{Application}
%
First we state our framework (\cite{OS2}). 
Let $M$ be an $m$-dimensional compact orientable submanifold of $\mathbb R^n$ and $\lambda$ be a real number. 
Put 
\[\Delta_\varepsilon=\{(x,y)\in \mathbb R^n\times \mathbb R^n\,|\,|x-y|<\varepsilon\}.\]
It is interesting to see the asymptotics of 
\[
E_{r^\lambda}(\varepsilon,M)=\iint_{M\times M\setminus\Delta_\varepsilon}|x-y|^\lambda\,d^mxd^my,
\]
where $d^mx$ and $d^my$ denote the standard Lebesgue measure of $M$. 
To be precise, we expand the above in a series in $\varepsilon$ and study the coefficients. 
For example, if $K$ is a knot in $\mathbb R^3$ and $\lambda=-2$ then 
\[
E_{r^{-2}}(\varepsilon,K)=\frac{2L(K)}{\varepsilon}+E(K)+O(\varepsilon),
\]
where $L(K)$ is the length of the knot and $E(K)$ denotes the knot energy given by \eqref{def_V_E}, and if $S$ is a closed surface in $\mathbb R^3$ and $\lambda=-4$ then 
\[
E_{r^{-4}}(\varepsilon,S)=\frac\pi{\varepsilon^2}\,A(S)-\frac\pi8\,\log\varepsilon\int_S\Delta(x)\,d^2x+E(S)
-\frac\pi{16}\int_S\Delta(x)\log\Delta(x)\,d^2x
-\frac{\pi^2}2\chi(S),
\]
where $\chi(S)$ is the Euler characteristic of $S$, and if $\Omega$ is a $2$-dimensional compact submanifold of $\mathbb R^2$ then 
\[
E_{r^{-4}}(\varepsilon,\Omega)=\frac{\pi}{\varepsilon^2}A(\Omega)-\frac{2}{\varepsilon}L(\partial\Omega)+E_{OS}(\Omega)-\frac{\pi^2}{4}\chi(\Omega)+O(\varepsilon), 
\]
where $A(\Omega)$ is the area of $\Omega$ and $E_{OS}$ is the energy defined in \cite{OS}. We remark that $E_{OS}$ is also invariant under M\"obius transformations. 
We conjecture that a similar formula holds for compact bodies in $\mathbb R^3$. 

Let us focus on the constant term of the series of $E_{r^\lambda}(\varepsilon,M)$, which, after some modification if necessary, we call the {\em renormalized $r^\lambda$-potential energy} of $M$, denoted by $E_{r^\lambda}(M)$. 

Now we can define functionals for knots as follows. 
Let $N_\varepsilon(K)$ be an $\varepsilon$-tubular neighbourhood of $K$. 
Expand $E_{r^\lambda}(\partial N_\varepsilon(K))$ and $E_{r^{\lambda'}}(N_\varepsilon(K))$ in series of $\varepsilon$. 
We conjecture that functionals that can capture global properties of knots appear as a coefficient of $\varepsilon^2$-term of $E_{r^\lambda}(\partial N_\varepsilon(K))$ 
and as a coefficient of $\varepsilon^4$-term of $E_{r^{\lambda'}}(N_\varepsilon(K))$. 

When $K_\circ$ is a round circle with radius $1$, our main theorem implies 
\begin{eqnarray*}
E_{r^{-4}}(\partial N_\varepsilon(K_\circ))
&=&\frac{\pi^3\varepsilon}{2\sqrt{1-\varepsilon^2}} \left(\frac{3\log 2 -1}{\varepsilon^2}  +2-{2\varepsilon^2}\right)\\
&=&\frac{\pi^3(3\log2-1)}{2\varepsilon}+\frac{3\pi^3\left(\log2+1\right)}{4}\varepsilon
+\frac{\pi^3\left(9\log2-11\right)}{16}\varepsilon^3+O(\varepsilon^5). 
\end{eqnarray*}
Therefore, the functional thus obtained from the $\varepsilon^2$-term of $E_{r^{-4}}(\partial N_\varepsilon(K))$ vanishes for round circles.

\end{document}